\documentclass[reqno]{amsart}
\usepackage{amssymb}
\usepackage{hyperref}
\usepackage{graphicx}
\usepackage{booktabs}
\usepackage{multirow}

\newtheorem{theorem}{Theorem}

\newtheorem{corollary}[theorem]{Corollary}

\theoremstyle{remark}
\newtheorem{remark}[theorem]{Remark}

\numberwithin{equation}{section}

\begin{document}

\title[Approximation of Askey-Wilson roots]
{Approximation of Askey-Wilson roots}

\author{Jan Felipe van Diejen}
\author{Andrés Soledispa}
\author{Adri\'an Vidal}

\address{
Instituto de Matem\'atica, Universidad de Talca,
Casilla 747, Talca, Chile}

\email{diejen@utalca.cl, andres.soledispa@utalca.cl, adrian.vidal@utalca.cl}

\subjclass[2010]{Primary: 33D45; Secondary: 26C10, 33F05, 47J26}
\keywords{Askey-Wilson polynomials, roots of orthogonal polynomials, fixed-point iteration}

\thanks{This work was supported in part by the {\em Fondo Nacional de Desarrollo
Cient\'{\i}fico y Tecnol\'ogico (FONDECYT)} Grant  \# 1250427  and the \emph{Beca de Doctorado Nacional ANID} Scholarship \# 21241601.}

\date{January 2025}

\begin{abstract}
This note presents a fixed-point formula designed to approximate the roots of Askey-Wilson poynomials for small parameter values.
\end{abstract}

\maketitle



\section{Introduction}
The Askey-Wilson polynomial occupies the top box of Askey's hierarchical scheme of basic hypergeometric  orthogonal polynomials; all other classical orthogonal polynomials arise as its specializations through parameter degenerations and/or limit transitions  \cite{ask-wil:some,koe-les-swa:hypergeometric}. 
The roots of the Askey-Wilson polynomials were studied at an early stage in connection with the (in)existence of perfect codes in families of classical association schemes
\cite{chi:zeros}. On the other hand, for many classical hypergeometric orthogonal polynomials a well-established and fruitful technique to locate the positions of the roots is through Stieltjes' electrostatic interpretation
\cite{sti:sur,sze:orthogonal}, cf. e.g. \cite{for-rog:electrostatics,gru:variations,ism:electrostatics,ism:classical,mar-mar-mar:electrostatic,mar-ori-mar:electrostatic,saf-tot:logarithmic} for important examples of contemporary progress along these lines. In this same spirit---while also exploiting additional links with the Bethe-Ansatz method for quantum integrable systems \cite{gau:bethe,kor-bog-ize:quantum}---recently the positions of the Askey-Wilson roots were characterized in terms of the global minimum of a strictly convex Morse function \cite{die:gradient,die-ems:solutions}. The purpose of the present note is to point out that for small parameter values this critical point encoding the positions of the Askey-Wilson roots can be readily approximated by  means of  a fixed-point iteration formula. We will estimate the error of the approximation by refining previous bounds on the positions of the Askey-Wilson roots stemming from \cite{die-ems:solutions} and also briefly compare these refined bounds with alternative bounds for the extremal roots of the Askey-Wilson polynomial highlighted in \cite{ken-jor:characterization}.

The remainder of this introduction  recalls a minimal compendium of requisites from the literature permitting to present our results.
The Askey-Wilson polynomial  $\mathrm{p}_n(\cos\theta )=\mathrm{p}_n(\cos\theta ; a_1,a_2,a_3,a_4|q)$ \cite[Chapter 14.1]{koe-les-swa:hypergeometric}  is a polynomial of degree
$n$ in $\cos\theta$ solving the $q$-difference  eigenvalue equation
\begin{align}\label{DE}
&A(\theta) \bigl( \text{p}_n\bigl(\cos(\theta-\mathrm{i}\log q )\bigr) -\text{p}_n(\cos\theta )\bigr)+\\
&A(-\theta)  \bigl( \text{p}_n\bigl(\cos(\theta+\mathrm{i}\log q )\bigr) -\text{p}_n(\cos\theta) \bigr) \nonumber 
= E_n \text{p}_n(\cos\theta)
\end{align}
with
\begin{align*}
A(\theta )= \frac{\prod_{1\leq r\leq 4}(1-a_r e^{\mathrm{i}\theta})}{(1-e^{2i\theta})(1-q e^{2i\theta})} \quad\text{and}\quad
E_n=q^{-n}(1-q^n)(1- q^{n-1}{\textstyle \prod_{1\leq r\leq 4} }a_r)
\end{align*}
(as a rational identity in the parameters $a_r$ and $q$).  For $a_r $ in the open unit disc
\begin{equation}\label{U}
 \mathbb{U}=\{ z\in\mathbb{C}\mid |z|<1\} ,
\end{equation}
 with non-real parameters arising as complex conjugate pairs,
the Askey-Wilson polynomials extend smoothly in $q\in (-1,1)$ to a real-valued
 orthogonal basis  with respect to a bounded continuous positive weight function on the interval $0<\theta<\pi$ (cf. e.g. \cite[Eq. (14.1.2)]{koe-les-swa:hypergeometric}). The general theory of orthogonal polynomials  \cite[Chapter 3.3]{sze:orthogonal}  therefore guarantees that for this parameter regime
 the Askey-Wilson polynomial $\mathrm{p}_n(\cos\theta )$ has $n$ simple roots within the interval of orthogonality:
 \begin{subequations}
 \begin{equation}\label{factorization}
p_n(\cos(\theta))= \kappa_n \prod_{1\leq j\leq n} \bigl( \cos\theta-\cos\theta_{j,n}\bigr) ,
\end{equation}
where the overall constant $\kappa_n$ is determined by the chosen normalization and the angles parameterizing the roots have been ordered such that
\begin{equation}\label{roots}
0<\theta_{1,n}<\cdots <\theta_{j,n}<\theta_{j+1,n}<\cdots <\theta_{n,n}<\pi.
\end{equation}
\end{subequations}
Upon substituting the factorization \eqref{factorization} into the difference equation \eqref{DE}, it follows that the roots \eqref{roots} of $\mathrm{p}_n(\cos\theta )$ uniquely solve the following algebraic system of equations in the angles $0<\theta_1<\theta_2<\dots<\theta_n<\pi$ \cite[Theorem 2]{die:equilibrium} (cf. also \cite[Proposition 2.1]{bih-cal:properties},
\cite[Chapter 16.5]{ism:classical}, \cite[Section II]{oda-sas:equilibrium}, and \cite[Section 7]{wie-zab:algebraization}):
\begin{equation}
\prod_{1\leq r\leq 4} \frac{ (1-a_r e^{\mathrm{i}\theta_j})  } { ( e^{\mathrm{i}\theta_j}-a_r )   } 
  \prod_{\substack{1\le k \le n \\ k\neq j}}
       \frac{ (1-q e^{\mathrm{i}(\theta_j + \theta_k)}) (1-q e^{\mathrm{i}(\theta_j - \theta_k)})   }
             { ( e^{\mathrm{i}(\theta_j + \theta_k)}-q ) ( e^{\mathrm{i}(\theta_j- \theta_k)}-q )   } =1
\end{equation}
for $ j=1,\ldots ,n .$   After taking the logarithm on both sides and dividing by the imaginary unit  $\mathrm{i}$,  the latter system of equations is converted into the transcendental guise \cite{die-ems:solutions}:
\begin{subequations}
\begin{align}\label{bethe}
\boxed{\sum_{1\leq r\leq 4} v_{a_r}(\theta_j)  +
 \sum_{\substack{1\leq k\leq n \\ k \neq j}} \bigl( v_{q}(\theta_j+\theta_k )+v_{q} (\theta_j-\theta_k) \bigr) =2\pi j }
\end{align}
for $j=1,\ldots ,n$, where
\begin{equation}\label{vu}
\boxed{v_\epsilon(\theta):=\int_0^\theta u_\epsilon(\vartheta)\text{d}\vartheta\quad \text{with}\quad u_\epsilon(\vartheta):=\frac{1-\epsilon^2}{1-2\epsilon\cos(\vartheta)+\epsilon^2} 
\quad (\text{for}\ |\epsilon | <1).}
\end{equation}
\end{subequations}
It is instructive to emphasize at this point that the transcendental system \eqref{bethe}, \eqref{vu} amounts to the equations for the (unique) critical point given by the global minimum of the following strictly convex semi-bounded Morse function on $\mathbb{R}^n$ \cite{die-ems:solutions}:
\begin{align*}\label{VM}
V(\theta_1,\ldots ,\theta_n):= \sum_{1\leq j <k \leq n} \left( \int_0^{\theta_j+\theta_k} v_q (\vartheta) \text{d}\vartheta +  \int_0^{\theta_j-\theta_k} v_q (\vartheta ) \text{d}\vartheta\right) & \\
+ \sum_{1\leq j\leq n} \left(   \sum_{1\leq r\leq 4}  \int_0^{\theta_j} v_{a_r}(\vartheta)\text{d}\vartheta -2\pi j \theta_j \right) .&
  \nonumber
\end{align*}

The above summary of results from the literature  culminates into the following useful characterization of the Askey-Wilson roots (cf. \cite[Theorem 1]{die:gradient}). 

\begin{theorem}[\cite{die:gradient,die-ems:solutions}]\label{AWzeros:thm}
For $-1<q<1$ and $a_1,a_2,a_3,a_4$  taken inside the open unit disk $\mathbb{U}$ \eqref{U} with possible non-real parameters arising in complex conjugate pairs, one has that the $n$ ordered roots $ \theta_{j,n}$ \eqref{roots} of the Askey-Wilson polynomial $\mathrm{p}_n(\cos\theta ; a_1,a_2,a_3,a_4|q)$ provide the coordinates of the \emph{unique} solution for the transcendental system \eqref{bethe}, \eqref{vu}.
\end{theorem}

Theorem \ref{AWzeros:thm} will be the starting point for the present note.
Specifically, in Section \ref{sec21}  we employ the transcendental system in question to refine the bounds for the Askey-Wilson roots from \cite{die-ems:solutions}. Next, in Section \ref{sec22}, we present a fixed-point formula for the solution of the transcendental system and employ our refined bounds to estimate the error of approximation.
The proofs of the stated results are relegated to Sections \ref{sec31} and \ref{sec32}, respectively.

\section{Approximation of Askey-Wilson roots}

\subsection{Bounds for the Askey-Wilson Roots}\label{sec21}

The following theorem improves previous bounds for the roots of the Askey-Wilson polynomial from
\cite[Theorem 4.6]{die-ems:solutions}. In brief, the main issue with the former bounds is that these tend to become imprecise for roots corresponding to angles $\frac{\pi}{2}<\theta_{j,n}<\pi$. By treating the interval $(0,\pi)$ more symmetrically the bounds for such angles are improved significantly.

\begin{theorem}[Bounds for $ \theta_{j,n}$]\label{bounds:thm} Let $-1<q<1$ and let
the parameters $a_1,a_2,a_3,a_4$ be taken from the open unit disk $\mathbb{U}$ \eqref{U} with possible non-real parameters arising as complex conjugate pairs. If $\max\{ |a_1|,|a_2|,|a_3|,|a_4|,|q|)>0$, then the $n$ ordered roots $ \theta_{j,n}$ \eqref{roots} of the Askey-Wilson polynomial $\mathrm{p}_n(\cos\theta ; a_1,a_2,a_3,a_4|q)$ obey the following inequalities:
\begin{subequations}
\begin{equation}\label{aw-bounds:a}
 \theta_{j,n}^{(-)}  \leq \theta_{j,n} \leq  \theta_{j,n}^{(+)}  ,
\end{equation}
for  $j=1,\ldots ,n$, where
\begin{align}
 \theta_{j,n}^{(-)}&=
 \begin{cases}
 \pi j k_{+,n}^{-1} &\text{if}\ 1\leq j \leq j^{(-)}_{n},\\
 \pi  \bigl( 1-(n+1-j)k_{-,n}^{-1} \bigr)&\text{if}\ j^{(-)}_{n}< j\leq n,
 \end{cases}
 \\ \intertext{and}
 \theta_{j,n}^{(+)}&=
 \begin{cases}
 \pi j k_{-,n}^{-1}  &\text{if}\ 1\leq j \leq j^{(+)}_{n},\\
  \pi  \bigl(1-(n+1-j)k_{+,n}^{-1} \bigr) &\text{if}\  j^{(+)}_{n} < j\leq n,
 \end{cases}
\end{align}
with
\begin{equation}\label{aw-bounds:b}
 j^{(+)}_{n}=\frac{(k_{+,n}-n-1)k_{-,n}}{k_{+,n}-k_{-,n}},\quad
 j^{(-)}_{n}= \frac{k_{+,n}(n+1-k_{-,n})}{k_{+,n}-k_{-,n}} ,
\end{equation}
and
\begin{equation}\label{aw-bounds:c}
k_{\pm,n}=k_{\pm,n} (a_1,a_2,a_3,a_4 | q):=  (n-1)  \Bigl( {\textstyle \frac{1+|q|}{1- |q|}}\Bigr)^{\pm 1}+ 
\frac{1}{2} \sum_{1\leq r\leq 4}\Bigl({\textstyle \frac{1+|a_r|}{1- |a_r|}}\Bigr)^{\pm 1} .
\end{equation}
\end{subequations} 
\end{theorem}

Theorem \ref{bounds:thm} merges two type of bounds for the Askey-Wilson roots:
\begin{subequations}
\begin{equation}\label{ineq:a}
 \pi j  k_{+,n}^{-1}  \leq \theta_{j,n} \leq  \pi j k_{-,n}^{-1}
\end{equation}
and
\begin{equation}\label{ineq:b}
 \pi \bigl(1-(n+1-j)k_{-,n}^{-1}\bigr) \leq \theta_{j,n} \leq  \pi \bigl(1-(n+1-j)k_{+,n}^{-1}\bigr) .
\end{equation}
\end{subequations}
The first inequality \eqref{ineq:a} amounts to the previous inequality (4.17a) of \cite[Theorem 4.6]{die-ems:solutions} (cf. also
\cite[Remark 5.5]{die-ems:solutions}). 
To prove Theorem \ref{bounds:thm} it is needed  to verify the second inequality \eqref{ineq:b}; this verification is relegated to Section \ref{sec31} below,
together with a brief supplementary analysis to select the sharpest of the inequalities for each value of $j$.

As an illustration, Table \ref{AWbounds:table} displays the bounds of Theorem \ref{bounds:thm} in a small numerical example.
\emph{Here and below unboxed bounds stem from \eqref{ineq:a} and boxed bounds stem from \eqref{ineq:b}.}
 It is instructive to compare with \cite[Table 1]{die-ems:solutions}, where one finds a corresponding table
showing the original (unboxed) bounds of the form \eqref{ineq:a} for this same example. The improvements of the bounds stemming from Theorem \ref{bounds:thm} are thus encoded by the boxed  bounds for $\theta_{3,5}$,  $\theta_{4,5}$ and  $\theta_{5,5}$.

\begin{table}[hbt]
\centering
\caption{Roots and their bounds from Theorem \ref{bounds:thm} for the Askey-Wilson polynomial $\text{p}_5\bigl(\cos\theta ;\frac{3}{10},-\frac{1}{5},\frac{3}{20},\frac{1}{10}|\frac{1}{10}\bigr)$.}
\label{AWbounds:table}
\begin{tabular}{@{}|lc|ccccc|@{}}
\midrule 
                                     &    $n=5$  & $j=1$ & $j=2$ & $j=3$ & $j=4$&       $j=5$           \\ \midrule

                &$\theta_{j,n}^{(-)}$  & 0.3999                       & 0.7999                      & 1.1998                      & \boxed{1.7915}  & \boxed{2.4666}  \\ \cmidrule(r){1-2} 
                 &$\theta_{j,n}$ & 0.4959                      &  0.9967                      & 1.5078                      & 2.0332  & 2.5773    \\ \cmidrule(r){1-2}
 &$\theta_{j,n}^{(+)}$  & 0.6750                       & 1.3501                     & \boxed{1.9418 }                     & \boxed{2.3417}  & \boxed{2.7416} \\ 
   \midrule
  \end{tabular}
\end{table}

For the extreme roots of the Askey-Wilson polynomial alternative bounds have been reported in \cite[Theorem 4.3]{ken-jor:characterization}. Specifically, for the Askey-Wilson polynomial
$\text{p}_n(x ;a,b,c,d|q)$ the formulas due to Kenfack Nangho and Jordaan provide a lower bound for the largest root
$\cos\theta_{1,n}$ and an upper bound for the smallest root $\cos\theta_{n,n}$. 
For comparison, Table \ref{extreme-bounds:table}
shows the corresponding bounds stemming from Theorem \ref{bounds:thm} for the extreme roots pertaining  to the example displayed in \cite[Table 1]{ken-jor:characterization}. It is clear from the formulas for the bounds provided by Kenfack Nangho and Jordaan that for $n \to \infty$ the value of the pertinent bound stabilizes towards a (quite elementary) limiting value. In the present example the limiting values are reached within a precision of $4$ decimals already for $n=7$. Indeed, with this precision in all cases of Table
\ref{extreme-bounds:table} the corresponding inequalities from \cite[Theorem 4.3]{ken-jor:characterization} take the form 
$\cos\theta_{1,n}\geq 0.9488$ and $\cos\theta_{n,n}\leq 0.3370$, respectively. In conclusion, in this example our upper bound
for the smallest root is significantly better than the upper bound from  \cite[Theorem 4.3]{ken-jor:characterization}, whereas our lower bound for the largest root is better only for $n$ sufficiently large (as of $n>12$ to be precise in this case).

\begin{table}[hbt]
\centering
\caption{Largest/smallest root and its lower/upper bound from Theorem \ref{bounds:thm} for the Askey-Wilson polynomial $\text{p}_n\bigl(x ;\frac{6}{7},\frac{5}{7},\frac{4}{7},\frac{3}{7}|\frac{1}{9}\bigr)$.}
\label{extreme-bounds:table}
\begin{tabular}{@{}|lc|cccc|@{}}
\midrule 
                                     &     & $n=7$ & $n=9$ & $n=12$   & $n=20$        \\ \midrule
 &$ \cos\theta_{1,n}$ & 0.9819                      & 0.9861                     & 0.9900              &    0.9948       \\ 
                &$\cos\theta_{1,n}^{(+)}$  & 0.8268                       & 0.8969                      & 0.9430            &   0.9799        \\ \cmidrule(r){1-2} 
 &$\cos\theta_{n,n}^{(-)}$  & \boxed{-0.8268}                     &\boxed{-0.8969}                    & \boxed{-0.9430}       &     \boxed{-0.9799}          \\ 
 &$\cos\theta_{n,n}$  & -0.8643                      & -0.9225                     &    -0.9588    &   -0.9862        \\ 
   \midrule
  \end{tabular}
\end{table}

\begin{remark}\label{chebyshev:rem} For notational convenience,
Theorem \ref{bounds:thm} excludes the elementary case when all parameters $a_1,a_2,a_3,a_4$ and $q$ vanish.  The Askey-Wilson polynomials then reduce to Chebyshev polynomials of the second kind and the inequalities \eqref{ineq:a}, \eqref{ineq:b} become tight:  $\theta_{j,n}\to\frac{\pi j}{n+1}$ for $j=1,\ldots ,n$.
\end{remark}

\subsection{Fixed-point Iteration towards the Askey-Wilson Roots}\label{sec22}
For $\epsilon\in\mathbb{U}$ \eqref{U} the function 
\begin{subequations}
\begin{equation}\label{nu:a}
\nu_\epsilon(\theta):=  v_\epsilon(\theta)-\theta
\end{equation}
 (cf. Eq. \eqref{vu})  is a  smooth odd
function of   $\theta\in\mathbb{R}$  that is periodic  with period $2\pi$.   Its Fourier series reads
\begin{equation}\label{nu:b}
 \nu_\epsilon(\theta)=\sum_{k >0} {\textstyle \frac{ 2 \epsilon^k}{k} }\sin (k\theta) \quad \bigl(\theta\in\mathbb{R}\bigr)
 \end{equation}
and the values on the fundamental interval $-\pi<\theta<\pi$ are conveniently computed  in terms of elementary  functions as follows:
\begin{equation}\label{nu:c}
 \nu_\epsilon(\theta)=2 \arctan \left( {\textstyle  \frac{1+\epsilon}{1-\epsilon}\tan \bigl( \frac{\theta}{2}\bigr)}\right)-\theta
 \quad \big(\theta\in (-\pi,\pi)\bigr) .
\end{equation}
\end{subequations}

Starting from the initial condition (cf. Remark \ref{chebyshev:rem})
\begin{subequations}
\begin{equation}\label{fp0}
\theta_{j,n}^{(0)}:=\frac{\pi j}{n+1}\quad  (j=1,\ldots, n),
\end{equation} 
let us define
a sequence of points
\begin{equation}\label{theta-nl}
\theta^{(l)}_n=\bigl( \theta^{(l)}_{1,n},\ldots , \theta^{(l)}_{j,n},\ldots, \theta^{(l)}_{n,n}\bigr)\in \mathbb{R}^n ,\quad\, l\in\mathbb{N}_0=\{l\in\mathbb{Z}\mid l\geq 0\} ,
\end{equation}
by means of the following recurrence:
\begin{align}\label{fpl+1}
    \theta_{j,n}^{(l+1)}:=\theta_{j,n}^{(0)}&-\frac{1}{2(n+1)}
     \sum_{1\leq r\leq 4} \nu_{a_r} \bigl(\theta_{j,n}^{(l)}\bigr)\\
    & -\frac{1}{2(n+1)} \sum_{\substack{1\leq k\leq n \\ k\neq j}} \Bigl(\nu_{q} \bigl(\theta_{j,n}^{(l)}+\theta_{k,n}^{(l)}\bigr) +\nu_{q} \bigl(\theta_{j,n}^{(l)}-\theta_{k,n}^{(l)}\bigr) \Bigr) \nonumber
\end{align}
\end{subequations}
for $j=1,\ldots,n$.

Our goal is to approximate the vector
\begin{equation}\label{theta-n}
\theta_n:=\bigl( \theta_{1,n},\ldots , \theta_{j,n},\ldots, \theta_{n,n}\bigr) 
\end{equation}
of ordered Askey-Wilson roots $\theta_{j,n}$ \eqref{roots} by $\theta_n^{(l)}$ \eqref{theta-nl}. The following theorem provides a bound
for the error of this approximation.

\begin{theorem}[Error of Approximation]\label{fixed-point-error:thm}
 Let $-1<q<1$ and let
the parameters $a_1,a_2,a_3,a_4$ be taken from the open unit disk $\mathbb{U}$ \eqref{U} with possible non-real parameters arising as complex conjugate pairs. For any $l\in\mathbb{N}_0$, the following inequality holds:
\begin{subequations}
\begin{equation}
   \| \theta_{n}- \theta^{(l)}_{n}\|\leq \rho_n^l \| \theta_{n}^{(+)}-\theta_{n}^{(-)}\| ,
\end{equation}
where
$\theta_n^{(\pm)}:=\bigl( \theta_{1,n}^{(\pm)},\ldots , \theta_{j,n}^{(\pm)},\ldots, \theta_{n,n}^{(\pm)}\bigr) $,
\begin{equation}\label{rho}
    \rho_n= \rho_n(a_1,a_2,a_3,a_4|q) := \Bigl(\frac{n-1}{n+1} \Bigr) \frac{2|q|}{1- |q|}+ 
\frac{1}{n+1} \sum_{1\leq r\leq 4}\Bigl({ \frac{|a_r|}{1- |a_r|}}\Bigr) ,
\end{equation}
and  $\| \cdot\|$ refers to the standard Euclidean norm in $\mathbb{R}^n$.
\end{subequations}
\end{theorem}
The proof of Theorem \ref{fixed-point-error:thm} hinges on a study of the stability of the  fixed point through an estimate for the Lipschitz constant of the recurrence; the details are worked out in Section \ref{sec32} below. The following practical convergence criterion ensuring our recurrence is a contraction is immediate from the above theorem.

\begin{corollary}[Convergence]\label{convergence:cor} Upon assuming parameter values as specified in  Theorem \ref{fixed-point-error:thm}, the following limits hold for all $1\leq j\leq n$:
\begin{equation}\label{convergence}
    \lim_{l\to\infty} \theta_{j,n}^{(l)}=\theta_{j,n} ,
\end{equation}
\emph{provided} $   \rho_n(a_1,a_2,a_3,a_4|q)<1$.
\end{corollary}

We conclude from  Corollary \ref{convergence:cor} that---for parameters values within the (standard) orthogonality domain---the ordered roots $\theta_{j,n}$ \eqref{roots}  of the Askey-Wilson polynomial
$\mathrm{p}_n(\cos\theta ; a_1,a_2,a_3,a_4|q)$ can be effectively approximated  by means of the fixed-point iteration \eqref{fp0}--\eqref{fpl+1}, \emph{provided} the parameters obey the additional contraction constraint $  \rho_n(a_1,a_2,a_3,a_4|q)<1$.
In order for this convergence condition on the parameters to be satisfied for any $n\in\mathbb{N}$,
 it is \emph{necessary} that $|q|<\frac{1}{3}$ (because of the limit  $\lim_{n\to\infty} \rho_n(a_1,a_2,a_3,a_4|q)=\frac{2|q|}{1-|q|}$) while it is \emph{sufficient} that in addition $|a_r|<\frac{1}{3}$ for $1\leq r\leq 4$ (since $\rho_n(\frac{1}{3},
\frac{1}{3},\frac{1}{3},\frac{1}{3}|\frac{1}{3})=1$).  In practice, the convergence is thus guaranteed for $n$ sufficiently large provided $|q|<\frac{1}{3}$.
In Table \ref{approx:table} our fixed-point method is tested on the example from Table \ref{AWbounds:table} (in which case $\rho_n=\rho_5\bigl(\frac{3}{10},-\frac{1}{5},\frac{3}{20},\frac{1}{10}|\frac{1}{10}\bigr)=\frac{883}{2856} <1$).

\begin{table}[hbt]
\centering
\caption{Fixed-point approximations \eqref{fp0}--\eqref{fpl+1} for the roots of the Askey-Wilson polynomial $\text{p}_5\bigl(\cos\theta ;\frac{3}{10},-\frac{1}{5},\frac{3}{20},\frac{1}{10}|\frac{1}{10}\bigr)$.}
\label{approx:table}
\begin{tabular}{@{}|lc|ccccc|@{}}
\midrule 
                                     &    $n=5$  & $j=1$ & $j=2$ & $j=3$ & $j=4$&       $j=5$           \\ \midrule

                &$\theta_{j,n}^{(0)}$  & 0.5236                     & 1.0472                     & 1.5708                     & 2.0944  & 2.6180  \\ \cmidrule(r){1-2} 
                 &$\theta_{j,n}^{(1)}$ & 0.4976                      &  1.0005                     & 1.5137                     & 2.0401  & 2.5830    \\ \cmidrule(r){1-2}
 &$\theta_{j,n}^{(2)}$  & 0.4961                       & 0.9971                     & 1.5084                     & 2.0339  & 2.5780 \\ \cmidrule(r){1-2}
                 &$\theta_{j,n}$ & 0.4959                      &  0.9967                      & 1.5078                      & 2.0332  & 2.5773    \\
   \midrule
  \end{tabular}
\end{table}

\begin{remark}\label{as-for:rem}
For $l=0$ the fixed-point iteration \eqref{fp0}--\eqref{fpl+1} gives rise to a first approximation formula for the ordered Askey-Wilson roots
$\theta_{j,n}$ \eqref{roots} in closed form:
\begin{align}\label{asymptotic-roots}
    \theta_{j,n}^{(1)}=&{\textstyle \frac{\pi j}{n+1}}
    -{\textstyle \frac{1}{2(n+1)}}
    \left( \sum_{1\leq r\leq 4} \nu_{a_r} \bigl({\textstyle \frac{\pi j}{n+1}}\bigr) + \sum_{\substack{1\leq k\leq n \\ k\neq j}} \Bigl(\nu_{q} \bigl({\textstyle \frac{\pi (j+k)}{n+1}}\bigr) +\nu_{q} \bigl({\textstyle \frac{\pi (j-k)}{n+1}}\bigr) \Bigr)\right) ,
    \end{align}
together with a corresponding error estimate $   \| \theta_{n}- \theta^{(1)}_{n}\|\leq \rho_n \| \theta_{n}^{(+)}-\theta_{n}^{(-)}\| $ stemming from Theorem \ref{fixed-point-error:thm}. Table \ref{relative-error:table} tests the asymptotics of the relative error $\varepsilon_{j,n}:=|\theta_{j,n}^{(1)}-\theta_{j,n}|/\theta_{j,n}$ of this first approximation as the degree $n$ increases.
\end{remark}

\begin{table}[hbt]
\centering
\caption{Asymptotics of the relative error $\varepsilon_{j,n}=|\theta_{j,n}^{(1)}-\theta_{j,n}|/\theta_{j,n}$ of the approximation formula \eqref{asymptotic-roots} for selected roots of the
Askey-Wilson polynomial $\text{p}_n\bigl(\cos\theta ;\frac{3}{10},-\frac{1}{5},\frac{3}{20},\frac{1}{10}|\frac{1}{10}\bigr)$.}
\label{relative-error:table}
\begin{tabular}{@{}|lc|cccc|@{}}
\midrule 
                                     &    & $n=10$ & $n=20$ & $n=40$ & $n=80 $       \\ \midrule

                &$\varepsilon_{1,n}$ & $2.6\cdot 10^{-3} $                   & $1.6\cdot 10^{-3} $                      & $8.7\cdot 10^{-4} $                      & $4.6\cdot 10^{-4} $     \\ \cmidrule(r){1-2} 
                 &$\varepsilon_{\frac{n}{2},n}$& $2.3\cdot 10^{-3} $                       &  $1.2\cdot 10^{-3} $                   & $6.1\cdot 10^{-4} $                   & $3.1\cdot 10^{-4} $      \\ \cmidrule(r){1-2}
 &$\varepsilon_{n,n}$  & $5.3\cdot 10^{-4} $                       & $1.2\cdot 10^{-4} $                   & $2.8\cdot 10^{-5} $                    &$6.5\cdot 10^{-6} $   \\ 
   \midrule
  \end{tabular}
\end{table}

\section{Proofs}

\subsection{Proof of Theorem \ref{bounds:thm}}\label{sec31}
For $z\in\mathbb{C}$ let $\text{Arg}(z)\in (-\pi,\pi]$ denote the principal value of the argument of $z$ with the additional convention that $\text{Arg}(0):=0$.
From the following identity for $u_\epsilon (\vartheta)$ \eqref{vu} with $\vartheta\in\mathbb{R}$
\begin{equation*}
 \mathfrak{Re} \bigl(u_\epsilon(\vartheta)\bigr)=  \frac{1}{2}  \Bigl(  u_{|\epsilon|}\bigl(\vartheta+ \text{Arg}(\epsilon)\bigr) + u_{|\epsilon|}\bigl(\vartheta- \text{Arg}( \epsilon)\bigr)  \Bigr) ,
\end{equation*}
it is immediate that \label{ub}
\begin{equation}
\forall \vartheta\in\mathbb{R},\epsilon\in\mathbb{U}:\qquad \frac{1-|\epsilon|}{1+|\epsilon|}\leq \mathfrak{Re} \bigl( u_\epsilon (\vartheta)\bigr) \leq \frac{1+|\epsilon|}{1-|\epsilon|}. 
\end{equation}
Hence, for the parameter values specified in Theorem \ref{bounds:thm} and any $j\in\{1, \dots, n\}$ the following inequalities hold $\forall\vartheta\in\mathbb{R}$:
    \begin{equation*}
        2k_{-,n} \leq\sum_{1 \leq r \leq 4} u_{a_r}(\vartheta) + \sum_{\substack{1\leq k\leq n \\ k \neq j}}\bigl(u_{q}(\vartheta+\theta_{k,n})+u_{q}(\vartheta-\theta_{k,n}) \bigr)  \leq 2k_{+,n}.
    \end{equation*}
  
    Upon integrating the latter inequalities over  $\vartheta\in (0,\theta_{j,n})$, one readily finds that
    \begin{align*}
      2  \theta_{j,n} k_{-,n} \leq 
      \underbrace{\sum_{1\leq r\leq 4}v_{a_r}(\theta_{j,n})  + \sum_{\substack{1\leq k\leq n \\ k \neq j}} \bigl( v_{q}(\theta_{j,n}+\theta_{k,n} )+v_{q} (\theta_{j,n}-\theta_{k,n}) \bigr) }_{\stackrel{\text{Theorem}~\ref{AWzeros:thm}}{=}\ 2\pi j} \leq 2 \theta_{j ,n}k_{+,n},
    \end{align*}
    which entails the inequalities \eqref{ineq:a}.
  
  On the other hand, when integrating instead  over  $\vartheta\in (\theta_{j,n},\pi)$ one arrives similarly at
    \begin{equation*}
       2(\pi-  \theta_{j,n}) k_{-,n} \leq 2\pi (n+1-j)
     \leq 2(\pi- \theta_{j ,n}) k_{+,n},
    \end{equation*}
    therewith settling the inequalities \eqref{ineq:b}. (Recall at this point that for $\theta\in\mathbb{R}$ and  $\epsilon\in\mathbb{U}$: $v_\epsilon(\theta+\pi)=v_\epsilon (\theta)+\pi$ and $v_\epsilon(-\theta)=-v_\epsilon(\theta)$, so in particular
    $v_\epsilon (\pi)=\pi$ and $v_\epsilon(\pi+\theta)+v_\epsilon(\pi-\theta)=2\pi$.)

    Finally, to deduce which of the inequalities in Eqs. \eqref{ineq:a} and \eqref{ineq:b} are actually sharpest for a given $j\in \{1,\ldots, n\}$, it suffices to observe the following elementary equivalences:
    \begin{align*}
        \pi \bigl(1-(n+1-j)k_{-,n}^{-1}\bigr) \lesseqgtr \pi jk^{-1}_{+,n} \quad \Longleftrightarrow \quad j \lesseqgtr j_{n}^{(-)},
    \end{align*}
    and 
    \begin{align*}
        \pi jk^{-1}_{-,n} \lesseqgtr \pi \bigl( 1-(n+1-j)k_{+,n}^{-1} \bigr) \quad \Longleftrightarrow \quad j \lesseqgtr j_{n}^{(+)},
    \end{align*}
therewith completing the proof of the inequalities asserted in Theorem \ref{bounds:thm}.

\subsection{Proof of Theorem \ref{fixed-point-error:thm}}\label{sec32}
When writing the  transcendental equations \eqref{bethe}, \eqref{vu} in terms of $\nu_\epsilon(\theta)$ \eqref{nu:a}--\eqref{nu:c}, the system is readily recasted as a fixed-point equation:
\begin{subequations}
\begin{equation}
\theta_j=U_j(\theta_1, \dots, \theta_n)\quad\text{for}\ j=1,\ldots,n,
\end{equation}
with
\begin{align}
       U_j(\theta_1, \dots, \theta_n):= \frac{\pi j}{n+1} & - \frac{1}{2(n+1)}\sum_{1\leq r\leq 4} \nu_{a_r}(\theta_j) 
        \\ &
        - \frac{1}{2(n+1)} \sum_{\substack{1\leq k\leq n \\ k \neq j}} \bigl( \nu_{q}(\theta_j+\theta_k ) + \nu_{q}(\theta_j-\theta_k) \bigr) . \nonumber
\end{align}
\end{subequations}
The entries of the Jacobi matrix $DU(\theta_1, \dots, \theta_n)$ for the fixed-point map
 $U(\theta_1, \dots, \theta_n) := \bigl(U_1(\theta_1, \dots, \theta_n), \dots, U_n(\theta_1, \dots, \theta_n)\bigr)$ read explicitly:
\begin{align*}
&  DU_{jk}(\theta_1, \dots, \theta_n)    := \partial_{\theta_k} U_j(\theta_1, \dots, \theta_n) =   \\ 
&        \begin{cases}
            1 - \frac{1}{2(n+1)}{\displaystyle \sum_{1\leq r\leq 4}} u_{a_r}(\theta_j) - \frac{1}{2(n+1)}{\displaystyle \sum_{\substack{1\leq i\leq n \\ i \neq j}} }\bigl( u_{q}(\theta_j+\theta_i ) + u_{q}(\theta_j-\theta_i) \bigr) &\text{if}\ j = k,\\
                    \frac{1}{2(n+1)}u_{q}(\theta_j-\theta_k) - \frac{1}{2(n+1)}u_{q}(\theta_j+\theta_k) &\text{if}\  j \neq k.
        \end{cases}
\end{align*}
Since $u_\varepsilon (\vartheta)$ \eqref{vu} is an even function it follows that
the Jacobi matrix in question is symmetric, so its spectral radius can be determined by means of the
Rayleigh-Ritz ratio:
\begin{align*}
    \lambda_{\min} = \min_{\substack{(x_1, \dots, x_{n}) \in \mathbb{R}^{n} \\ x_1^2+\cdots+x_n^2\neq 0}} \frac{\sum_{1 \leq j, k \leq n} x_jx_k DU_{jk}(\theta_1, \dots, \theta_n)}{x_1^2 + \cdots + x_n^2}
\end{align*}
and
\begin{align*}
    \lambda_{\max} = \max_{\substack{(x_1, \dots, x_{n}) \in \mathbb{R}^{n} \\ x_1^2+\cdots+x_n^2\neq 0}} \frac{\sum_{1 \leq j, k \leq n} x_jx_k DU_{jk}(\theta_1, \dots, \theta_n)}{x_1^2 + \cdots + x_n^2}.
\end{align*}

Upon inserting the uniform bounds
\begin{align*}
    \forall \vartheta\in\mathbb{R},\epsilon\in\mathbb{U}:\qquad -\frac{2|\epsilon|}{1-|\epsilon|} \leq 1 -  \mathfrak{Re} \bigl(u_\epsilon(\vartheta)\bigr) \leq \frac{2|\epsilon|}{1+|\epsilon|}
\end{align*}
 (cf. Eq. \eqref{ub}) into the explicit formula for the corresponding quadratic form
\begin{align*}
     &   \sum_{1 \leq j, k \leq n} x_jx_k DU_{jk}(\theta_1, \dots, \theta_n) = \frac{1}{2(n+1)}\sum_{1\leq j \leq n}\sum_{1 \leq r \leq 4} \bigl(1 - u_{a_r}(\theta_j)\bigr) x_j^2 
        \\ &
        + \frac{1}{2(n+1)}\sum_{1 \leq j < k \leq n} \Biggl(\bigl(1-u_{q}(\theta_j+\theta_k)\bigr)(x_j+x_k)^2 + \bigl(1-u_{q}(\theta_j-\theta_k)\bigr)(x_j-x_k)^2\Biggr)
\end{align*}
and the subsequent expansion of the squares on the second line, one sees after factoring out a common factor $x_1^2+\cdots +x_n^2$ that for the parameter values specified in Theorem \ref{fixed-point-error:thm}:
\begin{equation*}
 \lambda_{\min}\geq   -\Bigl(\frac{n-1}{n+1}\Bigr)\frac{2|q|}{1-|q|} - \frac{1}{n+1}\sum_{1\leq r\leq 4}\Bigl({ \frac{|a_r|}{1-|a_r|}}\Bigr) 
\end{equation*}
and
\begin{equation*}
\lambda_{\max}\leq \Bigl(\frac{n-1}{n+1} \Bigr)\frac{2|q|}{1+|q|}+\frac{1}{n+1}\sum_{1\leq r\leq 4}\Bigl({ \frac{|a_r|}{1+|a_r|}}\Bigr) .
\end{equation*}

We thus conclude that the spectral norm  of our symmetric Jacobi matrix
\begin{equation*}
\|  DU(\theta_1, \dots, \theta_n)\|:=\max_{\substack{x\in \mathbb{R}^n\\ \|x\|=1}} \| DU(\theta_1, \dots, \theta_n) x\| 
\end{equation*}
is given by a spectral radius $\rho\bigl(  DU(\theta_1, \dots, \theta_n)\bigr)=\max\bigl( |\lambda_{\max}|, |\lambda_{\min}| \bigr)$  that remains uniformly  bounded
by $ \rho_n=\rho_n(a_1,a_2,a_3,a_4|q)$ \eqref{rho}:
\begin{align*}
    \|  DU(\theta_1, \dots, \theta_n)\| \leq \rho_n.
\end{align*}
Our fixed-point map is therefore a Lipschitz map with a Lipschitz constant  bounded by $\rho_n$:
\begin{align*}
\forall \theta,\tilde{\theta}\in\mathbb{R}^n:\quad    \| U(\theta)- U(\tilde\theta)\| \leq \rho_n   \| \theta - \tilde\theta \|
\end{align*}
(cf. e.g. \cite[Chapter III.2]{edw:advanced}).
Evaluation at the fixed point $\theta=\theta_n$ and at the initial condition $\tilde\theta=\theta_n^{(0)}$ then entails upon iteration that
\begin{subequations}
\begin{equation}\label{ba}
\forall l\in\mathbb{N}_0:\quad    \| \theta_{n} - \theta^{(l)}_{n}\|\leq \rho_n^l \| \theta_{n} - \theta^{(0)}_{n}\|.
  \end{equation}
Moreover, from the inequalities \eqref{aw-bounds:a} in Theorem \ref{bounds:thm} and the elementary inequalities
$    \theta_{j,n}^{(-)}  \leq \theta_{j,n}^{(0)} \leq  \theta_{j,n}^{(+)}$ for $j=1,\ldots, n$ it is immediate that
\begin{equation}\label{bb}
 \| \theta_{n} - \theta^{(0)}_{n}\|\leq  \| \theta_{n}^{(+)}-\theta_{n}^{(-)}\| .
\end{equation}
\end{subequations}
The asserted error bound of Theorem \ref{fixed-point-error:thm} now follows by combining the inequalities in Eqs. \eqref{ba} and \eqref{bb}.

\bibliographystyle{amsplain}

\begin{thebibliography}{000000000}




\bibitem[AW85]{ask-wil:some}  R. Askey and J. Wilson,
Some basic hypergeometric orthogonal polynomials that generalize Jacobi polynomials,
Mem. Amer. Math. Soc. {\bf 54} (1985), no. 319, 55 pp.

\bibitem[BC16]{bih-cal:properties}  O. Bihun and F.  Calogero, Properties of the zeros of the polynomials belonging to the $q$-Askey scheme,
J. Math. Anal. Appl. {\bf 433} (2016), 525--542.



\bibitem[C87]{chi:zeros}  L. Chihara,
On the zeros of the Askey-Wilson polynomials, with applications to coding theory, 
SIAM J. Math. Anal. {\bf 18} (1987), 191--207. 


\bibitem[D05]{die:equilibrium} J.F. van Diejen, On the equilibrium configuration of the $BC$-type Ruijsenaars-Schneider system,
J. Nonlinear Math. Phys. {\bf 12} (2005), suppl. 1, 689--696.

\bibitem[D19]{die:gradient} J.F. van Diejen, Gradient system for the roots of the Askey-Wilson polynomial,
Proc. Amer. Math. Soc. {\bf 147} (2019), 5239--5249.

\bibitem[DE19]{die-ems:solutions} J.F. van Diejen and E. Emsiz, Solutions of convex Bethe Ansatz equations and the zeros of (basic) hypergeometric orthogonal polynomials,
Lett. Math. Phys. {\bf 109 } (2019), 89--112.

\bibitem[E73]{edw:advanced} C.H. Edwards, Jr., \emph{Advanced Calculus of Several Variables}, Academic Press, New York, 1973.

\bibitem[FR86]{for-rog:electrostatics} P.J. Forrester and J.B.  Rogers,
Electrostatics and the zeros of the classical polynomials, 
SIAM J. Math. Anal. {\bf 17} (1986), 461--468. 

\bibitem[G14]{gau:bethe} M. Gaudin,
{\em The Bethe Wavefunction}, Cambridge University Press, Cambridge, 2014.

\bibitem[G98]{gru:variations} F.A. Gr\"unbaum,
Variations on a theme of Heine and Stieltjes: an electrostatic interpretation of the zeros of certain polynomials,
J. Comput. Appl. Math. {\bf 99} (1998), 189--194. 

\bibitem[I00]{ism:electrostatics} M.E.H. Ismail, 
An electrostatics model for zeros of general orthogonal polynomials,
Pacific J. Math. {\bf 193} (2000), 355--369. 

\bibitem[I05]{ism:classical} M.E.H.  Ismail, {\em Classical and Quantum Orthogonal Polynomials in One Variable}, Encyclopedia of Mathematics and its Applications, vol. 98, Cambridge University Press, Cambridge, 2005.

\bibitem[KJ19]{ken-jor:characterization} 
M. Kenfack Nangho and K. Jordaan, A characterization of Askey-Wilson polynomials, Proc. Amer. Math. Soc. {\bf 147} (2019),
2465--2480.

\bibitem[KLS10]{koe-les-swa:hypergeometric}  R. Koekoek, P.A. Lesky, and R. Swarttouw,
{\em Hypergeometric Orthogonal Polynomials and their q-Analogues},
Springer Monographs in Mathematics. Springer-Verlag, Berlin, 2010.

\bibitem[KBI93]{kor-bog-ize:quantum} V.E. Korepin, N.M. Bogoliubov, and A.G. Izergin,
{\em Quantum Inverse Scattering Method and Correlation Functions},
Cambridge University Press, Cambridge, 1993.


\bibitem[MMM07]{mar-mar-mar:electrostatic}
F. Marcell\'an, A. Mart\'{\i}nez-Finkelshtein, and P. Mart\'{\i}nez-Gonz\'alez, Electrostatic models for zeros of polynomials: old, new, and some open problems, J. Comput. Appl. Math. {\bf 207} (2007), 258--272. 

\bibitem[MOS23]{mar-ori-mar:electrostatic} A. Martínez-Finkelshtein, R. Orive, J. Sánchez-Lara, 
Electrostatic partners and zeros of orthogonal and multiple orthogonal polynomials, Constr. Approx.   {\bf 58} (2023), 271--342.

\bibitem[OS05]{oda-sas:equilibrium} S. Odake and R. Sasaki,
Equilibrium positions, shape invariance and Askey-Wilson polynomials, J. Math. Phys. {\bf 46} (2005), no. 6, 063513, 10 pp.

\bibitem[ST24]{saf-tot:logarithmic} E.B. Saff and V. Totik, 
\emph{Logarithmic Potentials with External Fields}, Second Edition, 
Springer Nature, Cham, 2024.

\bibitem[S1885]{sti:sur} T.J. Stieltjes, Sur certains polyn\^{o}mes
qui v\'{e}rifient une \'{e}quation diff\'{e}rentielle lin\'{e}aire du second
ordre et sur la theorie des fonctions de Lam\'{e}, Acta Math. {\bf 6} (1885), 321--326. 

\bibitem[S75]{sze:orthogonal} G. Szeg\"o, {\em Orthogonal Polynomials}, Fourth Edition, American Mathematical Society, Colloquium Publications, vol. XXIII, American Mathematical Society, Providence, R.I., 1975.

\bibitem[WZ95]{wie-zab:algebraization}  P.B. Wiegmann and A.V. Zabrodin,
Algebraization of difference eigenvalue equations related to $U_q(sl_2)$,
Nuclear Phys. B   {\bf 451}(3) (1995), 699--724.

 
\end{thebibliography}

\end{document}